\documentclass{article}
\usepackage{graphicx} 
\usepackage{amsmath}
\usepackage{enumerate}
\usepackage{amsfonts}
\usepackage{algpseudocode}
\usepackage{algorithm}
\usepackage{amsthm}
\usepackage{xcolor}
\usepackage{tikz}
\usepackage{tkz-tab}
\usepackage{caption}
\usepackage{subcaption}
\usepackage{hyperref}
\usepackage{comment}

\title{Sparsity greedoids and pebble game algorithms for posets}
\author{Signe Lundqvist, Tovohery Randrianarisoa, Klara Stokes \\ and Joannes Vermant}

\begin{document}
\maketitle
\newtheorem{theorem}{Theorem}
\newtheorem{proposition}[theorem]{Proposition}
\newtheorem{lemma}[theorem]{Lemma}
\newtheorem{corollary}[theorem]{Corollary}
\newtheorem{rem
ark}[theorem]{Remark}
\newtheorem{definition}[theorem]{Definition}
\newtheorem{example}[theorem]{Example}
\newcommand{\out}{\textup{out}} 
\newcommand{\spa}{\textup{span}} 
\newcommand{\Stab}{\textup{Stab}} 
\newcommand{\peb}{\textup{peb}} 

\newcommand{\AcceptEdge}[1]{\textbf{Accept-Edge}(#1)}
\newcommand{\MovePebble}[1]{\textbf{Move-Pebble}(#1)}

\begin{abstract}
\noindent We generalise a sparsity condition for hypergraphs and show a result relating sparseness of hypergraphs to the decomposition of 
a modified incidence graph into edge-disjoint forests. We also give new sparsity conditions for posets and define an algorithm of pebble game type for posets to test when these sparsity conditions hold. Furthermore, we prove that under natural conditions, the sparsity conditions define a greedoid.  
\end{abstract}

\section{Introduction}

A graph $G=(V,E)$ is called \emph{$(k,l)$-sparse} if
\[
|E'| \le k|V(E')| - l \qquad \text{for all nonempty } E' \subseteq E,
\]
where $V(E')=\{v \in V ~|~ \exists e\in E'\text{ with } v\in e\}$ and the graph $G=(V,E)$ is called \emph{$(k,l)$-tight} if it is $(k,l)$-sparse and satisfies $|E| = k|V| - l$. When $0\leq l < 2k$, $(k,l)$-sparsity defines a matroid, known as the \emph{$(k,l)$-sparsity matroid} \cite{LOREA1979103}. 
\noindent
For some particular examples of $k$ and $l$, $(k,l)$-sparsity characterises important properties of the graph:
\begin{itemize}
    \item The bases of the $(1,1)$-sparsity matroid are spanning trees.
    \item By a theorem of Nash--Williams \cite{Nash-Williams}, the bases of the $(k,k)$-sparsity matroid are exactly the edge sets that decompose into $k$ edge-disjoint spanning trees.
    \item The bases of the $(2,3)$-sparsity matroid are the edge sets of generically minimally rigid graphs in $\mathbb{R}^2$ \cite{Laman, Geiringer}.
\end{itemize}

There are several generalisations of $(k,l)$-sparsity matroids. Streinu and Theran defined the following sparsity condition for $r$-uniform hypergraphs \cite{STREINU20091944}: an $r$-uniform  hypergraph $\mathcal{H} =  (V, E)$ is $(k, l)$-sparse if, for any nonempty subset of hyperedges $E'\subseteq E$, one has that
\begin{equation*}
    |E'| \leq k |V(E')| - l, 
\end{equation*}
where again $V(E')=\{v \in V ~|~ \exists e\in E'\text{ with } v\in e\}$. This leads to a matroid whenever $l$ satisfies $0\leq l\leq kr -1$.

\textit{Graded sparsity matroids} are a different generalisation of $(k,l)$-sparsity matroids. Let $G=(V,E)$ be a graph, and suppose that $V_0,\dots,V_n$ are subsets of $V$ such that $V_0 \subseteq V_1 \subseteq \dots \subseteq V_{n-1} \subseteq V_n=V$. Let $k$ be an integer and suppose that $L=(l_0,\dots,l_n)$ are integers such that $l_0 \geq l_1 \geq \dots \geq l_n$. A graph is graded $(k,L)$-sparse if $G_i=(V_i,E(V_i))$ is $(k, l_i)$-sparse. Graded sparsity matroids were introduced by Lee, Streinu, and Theran \cite{lee2007graded}. One such graded sparsity matroid characterises a rigidity-theoretic problem, namely line-constrained bar-joint frameworks in the plane \cite{Streinu2010}.

Another generalisation, \textit{count matroids}, may be found in \cite[Section 13.5]{Frank2012}. This definition generalises various matroids from discrete geometry, such as the $2$-dimensional rigidity matroid, but also count matroids on bipartite graphs \cite{Plane_Matroid}. One is given a graph $G=(V, E)$, a function $m: V\rightarrow \mathbb{Z}_{>0}$, together with some $l\in \mathbb{Z}_{\geq 0}$ such that $l < m(v) +m(u)$ for every edge $uv\in E$. Then one defines  $F\subseteq E$ to be independent if for every $V'\subseteq V$, one has
\begin{equation*}
    |F(V')|\leq \sum_{v\in V'} m(v) - l,
\end{equation*}
where $F(V') = \{e \in F ~\vert \exists v\in V' \textup{ with } v\in e\}$.
These form the independent sets of a matroid with ground set $E$. In particular, for bipartite graphs $G=(V_1, V_2, E)$, one has a matroid with ground set $E(V_1 \cup V_2)$, defined by 
\begin{equation*}
    k: V\rightarrow \mathbb{Z}_{>0} : v\mapsto  \begin{cases}
        k_1 \text{ if } v\in V_1\\
     k_2 \text{ if } v\in V_2
    \end{cases}
\end{equation*}
and $l \in \mathbb{Z}_{\geq 0}$, whenever $l < k_1 + k_2$. A set of edges $F$ is independent if, for any subset $F'\subseteq F$, one has $|F'| \leq k_1 |V_1\cap V(F')| +k_2 |V_2\cap V(F')| - l$. 

In this paper, we define a sparsity condition that is similar to graded sparsity but which turns out to exhibit rather different behaviour. Let $V=V_1\cup\cdots\cup V_n$ be a partition of the vertex set, with edges allowed only between $V_i$ and $V_{i+1}$ for the different values of $i$. For sequences of positive integers $K=(k_1,\dots,k_n)$ and $L=(l_{1,2},\dots,l_{n-1,n})$, an edge set $E'\subseteq E$ is called \emph{$(K,L)$–sparse} if every $V'\subseteq V$ satisfies
\[
|E'(V')|\;\le\;\sum_{i=1}^n k_i\,|V'\cap V_i| \;-\; l(V'),\]
with 
\[
l(V')=\sup\bigl\{\,l_{i,i+1}: V'\cap V_i\neq\emptyset,\ V'\cap V_{i+1}\neq\emptyset\,\bigr\},
\]
and \emph{$(K,L)$-tight} if equality holds for $V'=V$. In general, the $(K,L)$-sparse edge-sets are not the independent sets of a matroid, unlike in the case of graded sparsity matroids.

However, we will show that $(K,L)$-sparsity defines a \textit{greedoid}, which is a combinatorial structure that generalises matroids. Greedoids were first introduced by Korte and Lovász \cite{Greedoids_1}. For more information about greedoids and their applications, see the book by Korte, Lovász, and Schrader \cite{Greedoids_2}. 
\begin{definition}
    A \textit{greedoid language} is a set of strings $\mathcal{L}$ such that the following properties hold:
    \begin{enumerate}
        \item If $\alpha \in \mathcal{L}$ and $\beta$ is an initial segment of $\alpha$, then $\alpha \in \mathcal{L}$.
        \item For any $\alpha, \beta \in \mathcal{L}$ with $\vert \alpha \vert > \vert \beta \vert$, there is some $x\in \alpha$ such that $\beta x\in \mathcal{L}$. Here $\vert \alpha \vert$ denotes the length of the word $\alpha$.
    \end{enumerate}
\end{definition}

Let $G=(V,E)$ be a graph such that the vertex set has a partition $V=V_1\cup \dots \cup V_n$ and all edges $e\in E$ have endpoints in $V_{i} \cup V_{i+1}$ for some $1 \leq i \leq n$. Let $K=(k_1,...,k_n)$ and $L=(l_{1,2},...,l_{n-1,n})$ be sequences of integers and let $\mathcal{L}_{(K,L)}$ be the set of strings on $E$ which satisfy the following three properties:
\begin{enumerate}
\item{For any $\alpha=e_1\dots e_m$, the set $\{e_1, \dots, e_m\}$ is $(K,L)$-sparse}.
\item{For any $\alpha=e_1\dots e_m$ and $1\leq i< j \leq n$, if $e_{k_1} \in E(V_i \cup V_{i+1})$ and $e_{k_2} \in  E(V_j \cup V_{j+1})$, then $k_1 < k_2$.}
\item{For any $\alpha=e_1\dots e_m$ and $i\in \{1,\dots n-1\}$ such that there is an $e_j$ with $e_{j}\in E(V_i\cup V_{i+1})$, the edge set $\{e_1, \dots, e_m \} \cap E(V_1\cup \dots \cup V_i)$ is maximally $(K,L)$-sparse on $V_1\cup \dots \cup V_i$. }
\end{enumerate}

\noindent
We prove the following theorem in Section \ref{sec:greedoid}.

\begin{theorem}\label{thm:is_greedoid=true}
Suppose that the parameters $$K=(k_1, \dots, k_n)\in  \mathbb{Z}_{>0}^{n}, L= (l_{1,2}, \dots, l_{n-1,n})\in  \mathbb{Z}_{\geq 0}^{n-1},$$ satisfy the following inequalities:
\begin{align*}
    k_i + k_{i+1} &> l_{i, i+1}, \\
    k_{i+1} + l_{i-1, i} &> l_{i, i+1},\\
     l_{n-1,n} \geq  \dots &\geq l_{1,2}.
\end{align*} 
Let $G=(V,E)$ be a graph such that the vertex set has a partition $V=V_1\cup \dots \cup V_n$ such that all edges $e\in E$ have endpoints in $V_{i} \cup V_{i+1}$ for some $1 \leq i \leq n$. Then $\mathcal{L}_{(K,L)}$ is a greedoid language.
\end{theorem}

Greedoids are precisely the combinatorial structures to which the greedy algorithm applies.  There are polynomial time algorithms to test whether a set of edges is independent in a count matroid. The pebble game algorithm, introduced by Jacobs and Hendrickson, is an algorithm for determining whether a graph has a $(2,3)$-tight spanning subgraph \cite{Pebble_original}. The pebble game algorithm takes as its input a graph $G$ and outputs a maximum size $(2,3)$-sparse subgraph. In particular, the algorithm can be used to decide whether the input graph has a $(2,3)$-tight spanning subgraph. 

Lee and Streinu analysed the pebble game algorithm and showed that it can recognise $(k,l)$-sparse multigraphs \cite{LEE20081425}. Streinu and Theran extended the algorithm to recognise $(k,l)$-sparse $r$-uniform hypergraphs \cite{STREINU20091944}, and Lee, Streinu, and Theran extended it even further to recognise graded sparse graphs \cite{lee2007graded}. 

Berg and Jord\'an developed a degree-constrained orientation algorithm able to determine independence in count matroids on graphs, inspired by the pebble game algorithm \cite{Berg-Jordan}, which also applies to more general count matroids \cite[Section 13.5]{Frank2012}. 

There are also other algorithms for recognising sparse graphs. For example, another algorithm that can extract a $(2,3)$-tight subgraph was introduced by Gabow and Westerman \cite{Gabow-Westerman}, and an algorithm for a different sparsity condition based on max-flow-min-cut was developed by Imai and Sugihara \cite{IMAI198579}. See \cite{madarasi2025efficient} for a comparison of the different methods for algorithms recognising maximum size $(k,l)$-sparse graphs. 

In this paper, we show that there is a pebble game algorithm that recognises $(K,L)$-sparsity under the condition that the inequalities from Theorem \ref{thm:is_greedoid=true} apply. This pebble game algorithm is similar to those that recognise independent sets of sparsity matroids, even though the $(K,L)$-sparse graphs do not define a matroid but a greedoid. Details of this pebble game algorithm are given in Section \ref{sec:Algorithm}.

\section{Sparsity for bipartite graphs}\label{sec:sparsityrk2}

As a warm-up, we will consider sparsity conditions for bipartite graphs. The sparsity condition for posets that we will define later is a direct generalisation of the condition for bipartite graphs. These sparsity conditions on bipartite graphs also appear in rigidity theory \cite{Plane_Matroid, discrete_matroids}. We will extend a result by Servatius that relates a sparsity condition to tree packings \cite{k-planematroids}. This means that the sparsity condition enforces a strong combinatorial structure. Throughout this article, all graphs are allowed to be multigraphs. 

\begin{definition}\label{def: bipartite def}
    Let $G=(V_1, V_2, E)$ be a bipartite graph, with a partition of the vertex set into sets $V_1$ and $V_2$, and let $k_1, k_2 \in \mathbb{Z}_{>0}, l\in \mathbb{Z}_{\geq 0}$. Then we say that $G$ is $(k_1, k_2, l)$-sparse if, for every nonempty subset of $E' \subseteq E$, it holds that
\begin{equation*}
     \vert E'\vert \leq k_1 \vert V(E') \cap V_1 \vert +  k_2 \vert V(E')\cap V_2 \vert - l.
\end{equation*}
    We say that $G$ is $(k_1, k_2, l)$-tight if it is $(k_1, k_2, l)$-sparse, and if additionally
\begin{equation*}
     \vert E\vert = k_1 \vert V(E) \cap V_1 \vert +  k_2 \vert V(E)\cap V_2 \vert - l.
\end{equation*}
\end{definition}

As mentioned in the introduction, this defines a count matroid whenever $l < k_1 +k_2$. 

\begin{example}
    The case where $k_1 = d,$ $k_2 = 1$ and $l=d$ leads to the sparsity condition given by
    \begin{equation*}
        |E'| \leq d |V_1 \cap V(E')| + |V_2 \cap V(E')| - d.
    \end{equation*}
    This sparsity condition can be used to characterise the dimension of the space of liftings of a generic scene, or, dually, to characterise parallel redrawings of a generic hyperplane arrangement \cite{discrete_matroids}, and is called the $d$-plane matroid.
\end{example}

One can also compare Definition \ref{def: bipartite def} to the definition of sparsity for $r$-uniform hypergraphs from \cite{STREINU20091944}. Recall that for any hypergraph $\mathcal{H}=(V,E)$, one can define the incidence graph $I(\mathcal{H}) = (X, I)$, which has vertex set $X = V\cup E$ and edge set $I=\{\{v, e\} \subseteq X ~ \vert~ v\in e\}$. For $\lambda \in \mathbb{Z}_{>0}$, we denote by $\lambda \cdot I(\mathcal{H})$ the multigraph $(X, I')$ with $I' := \biguplus_{i=1}^{\lambda} I$, i.e., $I'$ is the multi-set consisting of $\lambda$ copies of each edge of $I(\mathcal{H})$. The following proposition is easy to verify, primarily using the fact that any edge $e\in E$ induces exactly $r$ incidences. 
\begin{proposition}\label{comparison-Hyp-Inc}
    Given $r,k,l \in \mathbb{Z}_{>0}$ with $r\geq 2$ and $0 < l < k r$. For all natural numbers $\lambda\geq 1$ such that $\lambda r -1 +k >l $, an $r$-uniform hypergraph $\mathcal{H}=(V,E)$ is $(k, l)$-sparse if and only if the bipartite graph $\lambda \cdot I(\mathcal{H})$  is $(k, \lambda r-1, l)$-sparse.
    Moreover, $\mathcal{H}$ is $(k, l)$-tight if and only if $\lambda\cdot I(\mathcal{H})$ is $(k, \lambda r-1, l)$ tight.
\end{proposition}

Thus, one recovers $(k,l)$-sparsity for hypergraphs from Definition \ref{def: bipartite def} applied to the incidence graph. However, we remark that the matroidal properties do not transfer cleanly since, in one case, there is a matroid with a base set being the hyperedge set of $\mathcal{H}$, and in the other case, there is a matroid with a base set being the edges of $\lambda \cdot I(\mathcal{H})$. 

To end this section, we will relate
$(k_1,k_2,l)$-sparsity to tree packings. The special case of $(k,1,k)$-tight and $(1,k,k)$-tight graphs was shown by Servatius \cite{k-planematroids}. 

Given a bipartite graph $G=(V_1,V_2,E)$, define its associated $(k_1,k_2)$-fold butterfly graph $B$ as follows: the vertex set consists of $k_1$ vertices for each element of $V_1$ and $k_2$ vertices for each element of $V_2$. For each edge $vw \in E$ with $v\in V_1$ and $w\in V_2$, there are edges between all $k_1$ vertices representing $v$ and all $k_2$ vertices representing $w$. 

\begin{theorem}
    Let $k_1, k_2, l\in \mathbb{Z}_{>0}$ with $k_1+k_2 >l$. A bipartite graph $G=(V_1,V_2,E)$ is $(k_1,k_2,l)$-sparse if and only if adding any $l-1$ edges between $V_1$ and $V_2$ results in a bipartite graph with an associated $(k_1,k_2)$-fold butterfly graph $B=(V, F)$ that decomposes into $k_1k_2$ edge-disjoint forests.
    \label{arboricity}
\end{theorem}
\begin{proof}
Suppose first that any resulting $B$ can be written as the union of $k_1k_2$ edge-disjoint forests. Then any such $B$ is $(k_1k_2, k_1k_2)$-sparse by the Nash-Williams theorem. For any nonempty subsets $V_1' \subseteq V_1$ and $V_2'\subseteq V_2$, adding the $(l-1)$ edges within $V_1' \cup V_2'$, one has, by the sparsity of the resulting graph $B$, that $$k_1k_2| E(V_1' \cup V_2')|+ (l -1)k_1k_2 \leq k_1k_2(k_1| V_1'|) + k_1k_2(k_1| V_1'|) - k_1k_2,$$ from which the $(k_1, k_2, l)$-sparsity of $G$ follows immediately.

 Conversely, suppose that we have added any $l-1$ edges to $G$, and let $B$ be the resulting butterfly graph. Let $V'\subseteq V(B)$. Note that all edges belonging to $F(V')$ are between $V_{1, i}'$ and $V_{2, j}'$ for some copies of subsets of $V_1$ and $V_2$ respectively. Since $G$ is sparse, we can bound $|F(V_{1, i}'\cup  V_{2, j}')|$ by sparsity:
    \begin{equation*}
        |F(V_{1, i}'\cup  V_{2, j}')| \leq k_1 |V_{1, i}'| + k_2 |V_{2, j}'| - 1
    \end{equation*}
    Then we have
    \begin{align*}
        |F(V')| &= \sum_{i=1}^{k_1} \sum_{j=1}^{k_2} |F(V_{1, i}'\cup  V_{2, j}')|\\
        &\leq  \sum_{i=1}^{k_1} \sum_{j=1}^{k_2}\left( k_1 |V_{1, i}'| + k_2 |V_{2, j}'| - 1\right)\\
        &=  k_2 k_1 \sum_{i=1}^{k_1} |V_{1, i}'| + k_1 k_2 \sum_{j=1}^{k_2} |V_{2, j}'| - k_1 k_2\\
        &=  k_2 k_1 (|V'| - 1).
    \end{align*}
    The result then follows by the Nash-Williams theorem \cite{Nash-Williams}. 
\end{proof}

\section{$(K,L)$-sparsity}\label{sec:greedoid}
\begin{definition}\label{KL-sparse}
    Let $G=(V,E)$ be a (multi)graph such that $V$ has a partition into subsets $V_1, V_2, \dots, V_n$, where every edge $e\in E$ has endpoints in $V_i$  and $V_{i+1}$ for some $i\in \{1,\dots, n-1\}$.
    
       Let $K=(k_1,\dots, k_n)\in \mathbb{R}_{>0}^{n}$ and ${L=(l_{1,2}, \dots , l_{n-1, n}) \in \mathbb{R}_{\geq 0}^{n-1}}$. For any $V'\subseteq V$, we define
    \begin{equation*}
        l(V') = \underset{\substack{
       i: V_{i}\cap V'\neq \emptyset\\ 
       V_{i+1}\cap V'\neq \emptyset}}
       {\sup} \{l_{i,i+1}\}.
    \end{equation*}
    We say that $G$ is $(K,L)$-sparse if, for any subset $V' \subseteq V$, one has
    \begin{equation*}
       \vert E(V')\vert \leq \sum_{i=1}^{n} k_i\vert V' \cap V_i\vert -l(V').
    \end{equation*}
    We note that the use of $\sup$ makes $l(V')=-\infty$ whenever the defining set is empty.\\ 
    Furthermore, $G$ is $(K,L)$-tight if $G$ is $(K,L)$-sparse, and additionally
        \begin{equation*}
       \vert E\vert = \sum_{i=1}^{n} k_i\vert V \cap V_i\vert -\sup_{i=1,\dots, n-1} \{l_{i,i+1}\}.
    \end{equation*}
    We furthermore say that a subset of edges $E'\subseteq E$ is $(K, L)$-sparse or $(K, L)$-tight if the graph $(V, E')$ is $(K,L)$-sparse or $(K,L)$-tight, respectively.  We say that a subset of edges $E'\subseteq E$ is maximally sparse if $E'$ is $(K,L)$-sparse, and if for any edge $vw \in E \setminus E'$ such that $v,w \in V(E')$, $E'\cup \{vw\}$ is not $(K,L)$-sparse.
\end{definition}

We will mainly consider the case where the entries of $K$ and $L$ are integers, and we need the following inequalities to hold for the $(K,L)$-sparse graphs to be a family of graphs that is well-behaved. 
\begin{align}
    k_i + k_{i+1} &> l_{i, i+1} \label{conditionposet1 -k k l},\\
    k_{i+1} + l_{i-1, i} &> l_{i, i+1}\label{conditionposet2 - k l l},\\
     l_{1,2} \leq  \dots &\leq l_{n-1,n}. \label{conditionposet3-increasing}
\end{align}
The first inequality ensures that a single edge is always sparse.
We will show shortly that the second inequality and the third inequality ensure that sparsity is preserved whenever one adds isolated vertices to the graph, even when adding a new part $V_{n}$. Moreover, in Section \ref{sec:Algorithm}, we will see that these inequalities are necessary to prove the correctness of the algorithm. The following example shows that if the inequalities do not hold, then a sparse graph together with an isolated vertex may be non-sparse.

\begin{example}
Suppose that $k_1 = 1, k_2=3, k_3=1$, and $l_{1,2}=1, l_{2,3}=3$. Consider the following graph, where the label on the node denotes which set $V_i$ each node belongs to.
\begin{center}
\begin{tikzpicture}
\node[circle,draw] (A) at (-2,1) {$1$};

\node[circle,draw] (B) at (-2, 0) {$1$};
\node[circle,draw] (C) at (-2, -1) {$1$};
\node[circle,draw] (X) at (-2, -2) {$1$};
\node[circle,draw] (Y) at (-2,+2) {$1$};

\node[circle,draw] (D) at (0,1) {$2$};
\node[circle,draw] (E) at (0,-1) {$2$};
\node[circle,draw] (F) at (2, 0) {$3$};
\draw (A) to (D);
\draw (A) to (E);
\draw (B) to (D);
\draw (B) to (E);
\draw (C) to (D);
\draw (C) to (E);
\draw (X) to (D);
\draw (X) to (E);
\draw (Y) to (D);
\draw (Y) to (E);

\end{tikzpicture}
\end{center}
One easily checks that for any subset
    $V' \subseteq V_1 \cup V_2$
it holds that
\begin{equation*}
 \vert E(V') \vert \leq  \vert V' \cap V_1 \vert  + 3 \vert V' \cap V_2 \vert - 1.
\end{equation*}
Yet, taking $V= V_1\cup V_2 \cup V_3$, one finds $\vert E(V)\vert =10$; however, the desired bound is
\begin{equation*}
\vert V' \cap V_1 \vert  + 3 \vert V' \cap V_2 \vert +  \vert V' \cap V_3 \vert - 3 = 5 + 6 +1 - 3 =9,
\end{equation*}
so the graph is not sparse. Note that inequality \eqref{conditionposet2 - k l l} does not hold, since $k_{3}+l_{1,2} =1+1 < 3 = l_{2,3}$.
\end{example}

\begin{lemma}
    Let $K=(k_1, \dots, k_n)\in  \mathbb{Z}_{>0}^{n}$ and $L= (l_{1,2}, \dots, l_{n-1,n})\in  \mathbb{Z}_{\geq 0}^{n-1}$ be such that $K$ and $L$ satisfy inequalities \eqref{conditionposet1 -k k l}, \eqref{conditionposet2 - k l l}, and \eqref{conditionposet3-increasing}. Let $G=(V,E)$ be a (multi)graph such that $V$ has a partition into subsets $V_1, V_2, \dots, V_n$ such that for all $e\in E$ there is some $i$ such that $e \in E(V_i, V_{i+1})$.
    
         Suppose that all edges are supported on $V_1\cup \dots \cup V_{t}$ with $t<n$. Then $G$ is $(K, L)$-sparse if and only if $G'= (V_1\cup \dots \cup V_{t}, E)$ is $(K',L')$-sparse, where $K'= (k_1, \dots, k_t)$ and $L'= (l_{1,2}, \dots, l_{t-1, t})$.
\end{lemma}
\begin{proof}
`Only if' is clear. Thus, suppose that $G'$ is a $(K',L')$-sparse graph. It suffices to show that sparsity is preserved when adding a set $W$ of isolated vertices of $V_{t+1}, \dots, V_n$ to a set $V' \subseteq V_{1}\cup \dots \cup V_t$. Assume first $W\cup V' \subseteq V_1\cup \dots \cup V_{t+1}$, with $V' \cap V_{t-1} \neq \emptyset$, $ V' \cap V_{t}\neq \emptyset$, and $W \cap V_{t+1} \neq \emptyset$. Then inequality \eqref{conditionposet2 - k l l} and the sparsity of $G'$ imply that
\begin{align*}
   \vert E(W\cup V') \vert&= \vert E(V') \vert \leq \sum_{i=1}^{t} k_i \vert V'\cap V_i\vert - l_{t-1,t}\\
&< \sum_{i=1}^{t} k_i \vert V' \cap V_i\vert + k_{t+1} - l_{t, t+1} \leq \sum_{i=1}^{t+1} k_i \vert (W\cup V') \cap V_i\vert - l_{t, t+1}, 
\end{align*}
and hence, since inequality \eqref{conditionposet3-increasing} ensures that $l(V' \cup W)=l_{t, t+1}$ in the last expression, this subset meets the sparsity condition.
If $V' \cup W$ is a set with $W \cap V_{t+j} \neq \emptyset$, $W \cap V_{t+j+1}\neq \emptyset$ for some $j > 0$, then by taking the largest such $j$, we see that 
\begin{align*}
\vert E(W\cup V') \vert&= \vert E(V') \vert \leq \sum_{i=1}^{t} k_i \vert V '\cap V_i\vert - l(V')\\
&<  \sum_{i=1}^{t} k_i \vert V '\cap V_i\vert - l(V') + (k_{t+j} + k_{t+j+1} - l_{t+j, t+j+1}) \\
&\leq \sum_{i=1}^{n} k_i \vert (V'\cup W) \cap V_i\vert - l(V'\cup W)
\end{align*}
where we used inequality \eqref{conditionposet1 -k k l} to go from the first line to the second. For the case where $(V'\cup W) \cap V_{t} \neq \emptyset$, $W \cap V_{t+1}\neq \emptyset$, and $V_{t-1}\cap V'= \emptyset$, one similarly uses  \eqref{conditionposet1 -k k l}, and in any other case, it is trivial to show that $V'\cup W$ satisfies the sparsity condition. 
\end{proof}

We can say a bit more by analysing the proof of the above lemma. We note that all of the inequalities showing sparsity were strict when there was an isolated vertex. If one were to add an isolated vertex to $V_i$, where $i\leq t$, then one would also clearly not have a tight graph.  Hence, we have the following lemma. 

\begin{lemma}\label{no_isolated}
    Let $K=(k_1, \dots, k_n) \in \mathbb{Z}_{>0}^{n}$ and $L=(l_{1,2}, \dots, l_{n-1,n}) \in  \mathbb{Z}_{\geq 0}^{n-1}$ be such that \eqref{conditionposet1 -k k l}, \eqref{conditionposet2 - k l l}, and \eqref{conditionposet3-increasing} hold, and let $G=(V,E)$ be a $(K,L)$-tight graph. Then $G$ has no isolated vertices.
\end{lemma}

Definition \ref{KL-sparse} is particularly applicable to the underlying graphs of Hasse diagrams of graded posets, and we will consider a class of examples of this type. A graded poset is a partially ordered set $(\mathcal{P}, \leq)$ such that every maximal chain has the same finite length $n$. This allows one to define a rank function
\begin{equation*}
    r: \mathcal{P} \rightarrow \{0, \dots, n\}: x\mapsto r(x)
\end{equation*}
where $r(x)$ is the unique length of a maximal chain contained in $\{y\in \mathcal{P} ~ \vert ~ y\leq x\}$. 

Given a graded poset $(\mathcal{P}, \leq)$, the underlying graph of the Hasse diagram is the graph $G(\mathcal{P})=(V, E)$, where $V=\mathcal{P}$, and
\begin{align*}
    E&=\cup_{k=0}^{n}\{ \{x,y\} \subseteq \mathcal{P} ~|~ x\leq y, r(x)=k, r(y)=k+1 \}
\end{align*}
The vertex set of $G(\mathcal{P})$ can be partitioned into sets $V_1, \dots, V_{n+1}$, where $V_k$ consists of the elements of $\mathcal{P}$ of rank $k-1$.
By the definition of $G(\mathcal{P})$, all edges are between $V_j$ and $V_{j+1}$ for some $j\in \{1,\dots, n\}$.

If the underlying graph of the Hasse diagram of a poset $\mathcal{P}$ is $(K,L)$-sparse or $(K,L)$-tight, then we simply say that $\mathcal{P}$ is $(K,L)$-sparse or $(K,L)$-tight, respectively.

\begin{example}
We consider the graph resulting from the Hasse diagram of a cube. The graph has vertex set $V_1 \cup V_2 \cup V_3$, where $V_1$ represents the vertices of the cube, $V_2$ represents the edges of the cube, and $V_3$ represents the faces of the cube. The order relation is given by inclusion. For any face $z$, the set $R^z= \{x\in \mathcal{P}~\vert~ x < z\}$ induces the face lattice of a square. Hence, for $V'\subseteq R^z$, one has
\begin{equation*}
    \vert E(V')\vert \leq \vert V'\cap V_1 \vert   + \vert V'\cap V_2\vert .
\end{equation*}
Since every vertex is contained in $3$ faces, every edge is contained in $2$ faces, and every incident vertex-edge pair appears in $2$ faces, for any $V'\subseteq V_1\cup V_2$, one has
\begin{align*}
    2 \vert E(V')\vert &\leq \sum_{z\in Z} \vert V'\cap V_1\cap R^z \vert   + \vert V'\cap V_2 \cap R^z\vert\\
    &\leq 3 \vert V'\cap V_1 \vert  + 2 \vert V'\cap V_2 \vert,
\end{align*}
and equality holds when $V'=V_1\cup V_2.$

Making a similar reasoning for the sets $R_{h}=\{x\in \mathcal{P}~\vert ~h < x\}$, we get for $V'\subseteq V_2 \cup V_3$ that
\begin{align*}
    2 \vert E(V')\vert &\leq 2 \vert V'\cap V_2 \vert   + 4 \vert V'\cap V_3\vert,
\end{align*}
and equality holds when $V'=V_2\cup V_3.$

For general subsets $V'\subseteq V_1\cup V_2 \cup V_3$, we then see that
\begin{align*}
    \vert E(V')\vert &\leq \frac{3}{2}\vert V'\cap V_1 \vert + 2 \vert V'\cap V_2 \vert   + 2 \vert V'\cap V_3\vert,
\end{align*}
and equality holds when $V'=V.$

Thus, we see that the face lattice of a cube is $(K,L)$-tight with $K=(\frac{3}{2}, 2, 2)$, and $L=(0,0)$. 
By the same argument, the face lattice of a regular convex $3$-polytope with Schläfli symbol $\{p,q\}$ is $(K, L)$-tight with parameters $K=(\frac{q}{2}, 2,\frac{p}{2})$ and $L=(0,0)$.
\end{example}

The following example shows that the $(K,L)$-sparse sets are not, in general, the independent sets of a matroid. 
\begin{example}
Consider $(K,L)$-sparsity with $K=(1,1,2)$ and $L=(1,2)$. The number on each vertex denotes which set $V_i$ the vertex belongs to.
\begin{center}    
Graph $1$:\\
\begin{tikzpicture}
\node[circle,draw] (A) at (0,1) {$1$};
\node[circle,draw] (B) at (0, 0) {$1$};
\node[circle,draw] (C) at (1, 0) {$2$};
\node[circle,draw] (D) at (1, 1) {$2$};
\node[circle,draw] (E) at (2, 0) {$3$};
\node[circle,draw] (F) at (2,1) {$3$};
\draw (A) to (D);
\draw (B) to (C);
\draw (C) to (E);
\draw (D) to (E);
\draw (C) to (F);
\draw (D) to (F);
\end{tikzpicture}\\
Graph $2$:\\
    \begin{tikzpicture}
\node[circle,draw] (A) at (0,1) {$1$};
\node[circle,draw] (B) at (0, 0) {$1$};
\node[circle,draw] (C) at (1, 0) {$2$};
\node[circle,draw] (D) at (1, 1) {$2$};
\node[circle,draw] (E) at (2, 0) {$3$};
\node[circle,draw] (F) at (2,1) {$3$};
\draw (A) to (D);
\draw (B) to (C);
\draw (B) to (D);
\draw (D) to (E);
\draw (D) to (F);
\end{tikzpicture}
\end{center}

Both graphs are $(K,L)$-sparse; however, one cannot add any edge from graph 1 to graph 2 without breaking sparsity. This means that the augmentation property for a matroid is not satisfied.
\end{example}

As stated in the introduction, we will show that one does, however, get a greedoid. We will consider strings $e_1\dots e_m$ on the alphabet $E$. For any string $\alpha=e_1\dots e_m$, we define the edge set $E_\alpha:= \{e_1, \dots, e_m\}$. Then, we define the collection $\mathcal{L}_{(K,L)}$ to be the set of strings on $E$ that satisfy the following three properties:
\begin{enumerate}
\item{For any $\alpha=e_1\dots e_m$, the set $E_\alpha$ is $(K,L)$-sparse}.
\item{For any $\alpha=e_1\dots e_m$ and $1\leq i< j \leq n-1$, if $e_{k_1} \in E(V_i \cup V_{i+1})$ and $e_{k_2} \in  E(V_j \cup V_{j+1})$, then $k_1 < k_2$.}
\item{For any $\alpha=e_1\dots e_m$ and $i\in \{1,\dots, n-1\}$ such that there is an $e_j$ with $e_{j}\in E(V_i\cup V_{i+1})$, the subgraph $E_{\alpha} \cap E(V_1\cup \dots \cup V_i)$ is maximally $(K,L)$-sparse on $V_1\cup \dots \cup V_i$. }
\end{enumerate}

We will use the terminology \textit{maximal feasible set} to mean a set $E_{\alpha}$, where $\alpha \in \mathcal{L}_{(K,L)}$ is a maximal word, meaning there is no word $\beta \in \mathcal{L}_{(K,L)}$ such that $\beta = \alpha e$ for some $e$. Expanding on what this means, the maximal feasible sets $E'$ are the sets $E'\subseteq E$ such that $E'(V_1 \cup \dots \cup V_j)$ is a maximally sparse set for each $j\in \{2, \dots, n\}$.

\begin{example}
    Consider the graph pictured below, where $V_1=\{v_1,v_2,v_3, v_4\}$, $V_2=\{x_1,x_2, x_3, x_4\}$ and $V_3=\{y_1,y_2,y_3,y_4\}$.
    \begin{figure}[h!]
    \centering
        \begin{tikzpicture}
\node[circle,draw] (A1) at (0,1) {$v_1$};
\node[circle,draw] (B1) at (0, 0) {$v_2$};
\node[circle,draw] (C1) at (0, -1) {$v_3$};
\node[circle,draw] (D1) at (0, -2) {$v_4$};
\node[circle,draw] (A2) at (2,1) {$x_1$};
\node[circle,draw] (B2) at (2, 0) {$x_2$};
\node[circle,draw] (C2) at (2, -1) {$x_3$};
\node[circle,draw] (D2) at (2, -2) {$x_4$};
\node[circle,draw] (A3) at (4,1) {$y_1$};
\node[circle,draw] (B3) at (4, 0) {$y_2$};
\node[circle,draw] (C3) at (4, -1) {$y_3$};
\node[circle,draw] (D3) at (4, -2) {$y_4$};
\draw (A1) to (A2);
\draw (A1) to (B2);
\draw (A1) to (C2);

\draw (B1) to (A2);
\draw (B1) to (B2);
\draw (B1) to (C2);
\draw (C1) to (C2);
\draw (C1) to (A2);
\draw (C1) to (B2);
\draw (D1) to (C2);

\draw (D1) to (D2);
\draw (A2) to (A3);
\draw (A2) to (B3);
\draw (B2) to (A3);
\draw (B2) to (B3);
\draw (B2) to (C3);
\draw (C2) to (B3);
\draw (C2) to (C3);
\draw (C2) to (D3);
\draw (D2) to (C3);
\draw (D2) to (D3);
\end{tikzpicture}
\end{figure}

We will consider what the maximal feasible sets are for $K=(1,2,2)$ and $L=(1,2)$. We thus want to build a subgraph such that it is maximally sparse on $V_1 \cup V_2$, and then add the maximum possible number of edges between $V_2$ and $V_3$. The condition 
\begin{equation*}
    |E(V')| \leq |V'\cap V_1 | + 2|V'\cap V_2| - 1
\end{equation*}
means we can pick at most $8$ edges supported on the set $V'=\{v_1,v_2,v_3,x_1,x_2,x_3\}$. Having picked $8$ edges $E_* = \{e_1, \dots, e_8\}$ supported on $V'$, the condition for sets $V" = V' \cup \{y_i\}$ gives
\begin{equation*}
    |E(V")| = |E_* \cup E"|  \leq |V"\cap V_1 | + 2|V"\cap V_2| + 2|V"\cap V_3| - 2
\end{equation*}
implies
\begin{equation*}
    |E"| \leq 2|V"\cap V_3| - 1,
\end{equation*}
where $E"$ is the set of edges in $E(V")$ supported on $V_2\cup V_3$. Hence, any $y_i$ can be connected to at most one of $x_1, x_2, x_3$. Analysing the other edges, we see that all of them can be added without breaking sparsity. Below, there are two examples of maximal feasible sets.
    \begin{figure}[h!]
    \centering
        \begin{tikzpicture}
\node[circle,draw] (A1) at (0,1) {$v_1$};
\node[circle,draw] (B1) at (0, 0) {$v_2$};
\node[circle,draw] (C1) at (0, -1) {$v_3$};
\node[circle,draw] (D1) at (0, -2) {$v_4$};
\node[circle,draw] (A2) at (2,1) {$x_1$};
\node[circle,draw] (B2) at (2, 0) {$x_2$};
\node[circle,draw] (C2) at (2, -1) {$x_3$};
\node[circle,draw] (D2) at (2, -2) {$x_4$};
\node[circle,draw] (A3) at (4,1) {$y_1$};
\node[circle,draw] (B3) at (4, 0) {$y_2$};
\node[circle,draw] (C3) at (4, -1) {$y_3$};
\node[circle,draw] (D3) at (4, -2) {$y_4$};

\draw (A1) to (B2);
\draw (A1) to (C2);
\draw (B1) to (A2);
\draw (B1) to (B2);
\draw (B1) to (C2);
\draw (C1) to (C2);
\draw (C1) to (A2);
\draw (C1) to (B2);
\draw (D1) to (C2);
\draw (D1) to (D2);

\draw (A2) to (A3);
\draw (A2) to (B3);
\draw (C2) to (C3);
\draw (C2) to (D3);
\draw (D2) to (C3);
\draw (D2) to (D3);

\node[circle,draw] (A1') at (6,1) {$v_1$};
\node[circle,draw] (B1') at (6, 0) {$v_2$};
\node[circle,draw] (C1') at (6, -1) {$v_3$};
\node[circle,draw] (D1') at (6, -2) {$v_4$};
\node[circle,draw] (A2') at (8,1) {$x_1$};
\node[circle,draw] (B2') at (8, 0) {$x_2$};
\node[circle,draw] (C2') at (8, -1) {$x_3$};
\node[circle,draw] (D2') at (8, -2) {$x_4$};
\node[circle,draw] (A3') at (10,1) {$y_1$};
\node[circle,draw] (B3') at (10, 0) {$y_2$};
\node[circle,draw] (C3') at (10, -1) {$y_3$};
\node[circle,draw] (D3') at (10, -2) {$y_4$};
\draw (A1') to (A2');
\draw (A1') to (B2');
\draw (A1') to (C2');
\draw (B1') to (A2');
\draw (B1') to (C2');
\draw (C1') to (C2');
\draw (C1') to (A2');
\draw (C1') to (B2');
\draw (D1') to (C2');
\draw (D1') to (D2');
\draw (B2') to (A3');
\draw (B2') to (B3');
\draw (B2') to (C3');
\draw (C2') to (D3');

\draw (D2') to (C3');
\draw (D2') to (D3');

\end{tikzpicture}
\end{figure}
\end{example}

To prove that $\mathcal{L}_{(K,L)}$ is a greedoid, we will show that the sets of edges that can be added to a word $\alpha$ between $V_j$ and $V_{j+1}$ are the independent sets of a matroid. Let us make this precise. Suppose that $E^*$ is a maximally sparse set supported on $V_1 \cup...\cup V_{j}$, and let $\mathcal{I}_j$ be the collection of subsets $E' \subseteq E(V_{j},V_{j+1})$ such that $E^* \cup E'$ is $(K,L)$-sparse. We will show that the elements of $\mathcal{I}_j$ are the independent sets of a matroid with ground set $E(V_{j},V_{j+1})$. We will moreover show (see Lemma \ref{Matroid_1}) that $\mathcal{I}_j$ does not depend on the chosen maximally sparse set. We denote the pair $(E(V_{j},V_{j+1}), \mathcal{I}_j)$ by $\mathcal{M}_j.$


\begin{theorem}\label{thm:Matroid}
Let $K  \in \mathbb{Z}_{>0}^{n}, L \in \mathbb{Z}_{\geq 0}^{n-1}$ satisfy inequalities \eqref{conditionposet1 -k k l}, \eqref{conditionposet2 - k l l}, and \eqref{conditionposet3-increasing}. Let $G=(V,E)$ be a graph such that the vertex set has a partition $V=V_1\cup \dots \cup V_n$ so that all edges $e\in E$ have endpoints in $V_{i} \cup V_{i+1}$ for some $1 \leq i \leq n$. Then $\mathcal{M}_{j}$ is a matroid for all $1\leq j \leq n-1$.
\end{theorem}

 Note that $\mathcal{M}_{1}$ is the count matroid with ground set $E(V_1\cup V_2)$, defined by 
\begin{equation*}
    k: V\rightarrow \mathbb{Z}_{>0} : v\mapsto  \begin{cases}
        k_1 \text{ if } v\in V_1\\
     k_2 \text{ if } v\in V_2
    \end{cases}
\end{equation*}
and $l \in \mathbb{Z}_{\geq 0}$. This will serve as the base case, and the proof of Theorem \ref{thm:Matroid} will proceed by induction. We will need several lemmas to prove Theorem \ref{thm:Matroid}.



We now introduce a key tool used to prove Theorem \ref{thm:Matroid}.  To simplify notation,
for any subset of edges $F\subseteq E$, we define $V_i(F):=V(F)\cap V_i$.
Let $\alpha_j$ be a word in $\mathcal{L}_{(K,L)}$ such that no $e\in E(V_{j+1} \cup V_{j+2})$ appears in $\alpha_j$, and $E_{\alpha_j}$ is maximally sparse on $V_1 \cup \dots \cup V_{j+1}$. We then define, for subsets $X\subseteq V_{j+1}$ the \textit{surplus}
\begin{align*}
    \theta_{\alpha_j}(X) := \min_{W: X\subseteq W}\{(\sum_{i=1}^{j+1} k_{i}|V_{i} \cap W| - |E_{\alpha_j}(W)|\}.
\end{align*} 
We remark that $\theta$ depends on the chosen parameters $K$.We will sometimes suppress the subscript for $\alpha_j$, writing $\alpha$ in stead of $\alpha_j.$ It follows easily from the definition that if $X\subseteq Y$ then $\theta_{\alpha}(X) \leq \theta_{\alpha}(Y)$.

Let us briefly motivate the definition of the surplus.
Fix some $\alpha \in \mathcal{L}_{(K,L)}$ such that $E_{\alpha} \subseteq E(V_1 \cup \dots \cup V_j)$ is a maximally sparse subset of $ E(V_1\cup \dots \cup V_j)$. Defining $\mathcal{I}_j$ using $E_{\alpha}$, one can see that a set $E'\subseteq E(V_j\cup V_{j+1})$ belongs to $\mathcal{I}_j$ if and only if
\begin{equation}
    |F|\leq k_{j+1}|V_{j+1}(F)| + \theta_{\alpha}(V_{j}(F)) - l_{j,j+1}
    \label{def:matroid}
\end{equation}
for any subset $F \subseteq E'$.  

For some classical count matroids, the surplus has interpretations in terms of fundamental graph invariants.

\begin{example}
Let $K$ be such that $k_1 = k_2 =1$ and $l_{1,2} = 1$.  Let $G=(V, E)$ be a $(K,L)$-sparse graph. This means that $\mathcal{M}_1$ is the graphic matroid.

Suppose that $\alpha=e_1 \dots e_m$ such that $e_i \in E(V_1 \cup V_2)$ for all $i\in \{1, \dots, m\}$, and that $E_{\alpha}$ is a basis of the graphic matroid on $V_1\cup V_2$, i.e., $T= (V_1\cup V_2, E_{\alpha})$ is a spanning forest. Then for $X\subseteq V_2$, one can verify that 
\begin{equation*}
\theta_{\alpha}(X) = | \{C \subseteq V_1\cup V_2 ~|~ C \textup{ is a connected component of } T \textup{ with } X\cap C \neq \emptyset \}|.  
\end{equation*}
In particular, if $T$ is a spanning tree, then $\theta_{\alpha}(X)= 1$ for all subsets $X\subseteq V_2$. One sees immediately that this definition is independent of the chosen $\alpha$.
\end{example}

\begin{example}
Let $K$ be such that $k_1 = k_2 =2$ and $l_{1,2} = 3$. 
Let $G=(V, E)$ be a $(K,L)$-sparse graph. This means that $\mathcal{M}_1$ is the $2$-dimensional rigidity matroid. Let $\alpha \in \mathcal{L}_{(K,L)}$ be such that $E_{\alpha}\subseteq E(V_1 \cup V_2)$ and such that $E_\alpha$ is a basis for $\mathcal{M}_1$. Then, $\theta_{\alpha}(X)$ gives the minimal degrees of freedom of subsets $W$ containing $X$. When $G$ is a rigid graph, one has $\theta_{\alpha}(X) = 3$ for all subsets with $|X|\geq 2$, and $\theta_{\alpha}(X) = 2$ for $X=\{v\}$. One sees that $\theta_{\alpha}$ is independent of the chosen $\alpha$.
\end{example}

Our first goal is to show that $\theta_{\alpha_j}$ is independent of the choice of $\alpha_j,$ whenever $\alpha_j$ is a maximal feasible set supported on $V_1\cup \dots \cup V_{j+1}$. The following lemma essentially states that the expression $\sum_{i=1}^{j+1} k_i|V_i \cap X| - |E_{\alpha_{j}}(X)|$ defining $\theta_{\alpha_j}$ decreases when adding a set of vertices $T$ such that $E_{\alpha_j}$ is tight on $T$ to the given set $X$.

\begin{lemma}\label{Increase_when_tight - greedoids}
Let $K\in \mathbb{Z}_{>0}^{n}, L\in \mathbb{Z}_{\geq 0}^{n-1}$ be such that inequalities \eqref{conditionposet1 -k k l}, \eqref{conditionposet2 - k l l} and \eqref{conditionposet3-increasing} are valid. Let $G=(V,E)$ be a graph such that the vertex set has a partition $V=V_1\cup \dots \cup V_n$ such that all edges $e\in E$ have endpoints in $V_{i} \cup V_{i+1}$ for some $1 \leq i < n$. 

Suppose also that $B \subseteq E$ is $(K,L)$-sparse. Let $X\subseteq V_{1}\cup \dots \cup V_{j}$ and $T \subseteq V_{1}\cup \dots \cup V_{j}$ be such that $X\cap T$ contains some vertices of $V_{j-1}$ and $V_{j}$. Suppose also that
\begin{equation*}
\sum_{i=1}^{j} k_i |V_i\cap T| - |B(T)|  = \ell_{j-1, j}. 
\end{equation*}
Then
\begin{equation*}
\sum_{i=1}^{j} k_i |V_i \cap X| - |B(X)|  \geq \sum_{i=1}^{j} k_i |V_i\cap (T\cup X)| - |B(T\cup X)|. 
\end{equation*}
\end{lemma}
\begin{proof}
Let $B(X, T)$ be the subset of edges of $B$ such that one endpoint is contained in $X$ and the other endpoint is contained in $T\setminus X$.
We have
\begin{align*}
&\hspace{11pt} \sum_{d=1}^{j} k_i |V_i\cap(T\cup X)| - |B(T\cup X)|\\
  &=  \sum_{d=1}^{j+1} k_i |V_i \cap (T\cup X)| - \left( |B(X)| + |B(T)| + |B(X,T)|- |B(X\cap T)|\right)\\
&=  \left( \sum_{d=1}^{j} k_i |V_i\cap X|  - |B(X)| \right) +  \left( \sum_{d=1}^{j+1} k_i |V_i\cap T|  - |B(T)| \right)\\
 & \hspace{11pt} - \left( \sum_{d=1}^{j} k_i |V_i \cap (X\cap T)|  - |B(X\cap T)| \right) - |B(X, T)|\\
\end{align*}

Since $B$ is $(K,L)$-sparse and since $X\cap T$ contains some vertices of $V_{j-1}$ and $V_{j}$, we have that $\sum_{d=1}^{j} k_i |V_i\cap (X\cap T)| -|B(X \cap T)| \geq l_{j-1,j}$. As $(T,B)$ is $(K,L)$-tight, it holds that:
\begin{align*}
\sum_{d=1}^{j} k_i |V_i\cap(T\cup X)| - |B(T\cup X)| &\leq  \sum_{d=1}^{j} k_i |V_i\cap X| - |B(X)| \\
& \vspace{10pt}+ \ell_{j-1,j} -\ell_{j-1,j} - |B(X, T)|\\
&\leq \sum_{d=1}^{j} k_i |V_i\cap X| - |B(X)| .
\end{align*}
This concludes the proof of the lemma.
\end{proof}

\begin{lemma}\label{blocks}
Let $K\in \mathbb{Z}_{>0}^{n}, L\in \mathbb{Z}_{\geq 0}^{n-1}$ be such that inequalities \eqref{conditionposet1 -k k l}, \eqref{conditionposet2 - k l l} and \eqref{conditionposet3-increasing} are valid. Let $G=(V,E)$ be a graph such that the vertex set has a partition $V=V_1\cup \dots \cup V_n$ such that all edges $e\in E$ have endpoints in $V_{i} \cup V_{i+1}$ for some $1 \leq i < n$, and assume that $G$ is $(K,L)$-sparse. 

Let $B_1, B_2 \subseteq E$ be tight edge sets, with $V(B_i)\subseteq V_1\cup \dots \cup V_t$. If $V(B_1)\cap V(B_2) \cap V_{t-1} \neq \emptyset$ and $V(B_1)\cap V(B_2) \cap V_{t} \neq \emptyset$, then $B_1 \cap B_2$ and $B_1\cup B_2$ are tight subgraphs.
Moreover, one has
\begin{align*}
    V(B_1\cup B_2) &= V(B_1) \cup V(B_2)\\
    V(B_1\cap B_2) &= V(B_1) \cap V(B_2)
\end{align*}
\end{lemma}

\begin{proof}
    Since sparsity holds by sparsity of $G$, we only need to show that $B_1\cup B_2$ and $B_1\cap B_2$ have the desired cardinality. We also remark that $V(B_1) \cup V(B_2) = V(B_1\cup B_2)$ is always true. By sparsity applied to $V(B_1) \cap V(B_2)$,
    \begin{equation*}
        \vert E(V(B_1) \cap V(B_2))\vert \leq \sum_{i=1}^{t} k_i \vert V(B_1)\cap V(B_2)\cap V_i \vert   -l_{t-1, t}.
    \end{equation*}
    Moreover, since any tight subgraph of a sparse graph is an induced subgraph, we see that $E(V(B_1) \cap V(B_2)) = B_1\cap B_2$. Indeed, if $e\in E(V(B_1) \cap V(B_2))$, then $e$ has endpoints in $V(B_1)$ and $V(B_2)$, and since $B_1$ and $B_2$ are induced, this means that if $e\in E(V(B_1) \cap V(B_2))$, then $e\in B_1\cap B_2$. 
    
    We then compute
    $$\begin{array}{cl}
     \vert B_1\cup B_2 \vert &= \vert B_1 \vert + \vert B_2 \vert -  \vert B_1\cap B_2 \vert\\
     &= \sum_{i=1}^{t} k_i\vert V_i(B_1) \vert  + \sum_{i=1}^{t} k_i\vert V_i(B_2) \vert  -2l_{t-1, t} - \vert B_1 \cap B_2 \vert \\
    &\geq \sum_{i=1}^{t} k_i \left( \vert V_i(B_1)\vert+ \vert V_i(B_2)  \vert \right)\\
     & ~~- \sum_{i=1}^{t} k_i \vert V_i(B_1)\cap V_i(B_2)\vert   -l_{t-1, t} \\
     &=\sum_{i=1}^{t} k_i\vert V_i(B_1\cup B_2) \vert  -l_{t-1, t},
\end{array}$$
    where we apply sparsity of $G$ to $B_1\cap B_2$ to get the third line from the second. By sparsity of $G$, it also holds that 
    \begin{equation*}
      \vert  B_1 \cup B_2 \vert \leq \sum_{i=1}^{t} k_i \vert V_i(B_1\cup B_2)\vert -l_{t-1, t},
    \end{equation*}
    and thus we have
    \begin{equation*}
    	  \vert  B_1 \cup B_2 \vert = \sum_{i=1}^{t} k_i \vert V_i(B_1\cup B_2)\vert  -l_{t-1, t}.
    \end{equation*} 
    Thus $B_1 \cup B_2$ is tight, and every inequality in the computation above is necessarily an equality. Thus, for $B_1 \cap B_2$, we see that
    $$\begin{array}{cll}
     \vert B_1\cap B_2 \vert & = & | B_1| +| B_2| - |B_1 \cup B_2| \\ 
    & = &\sum_{i=1}^{t} k_i \big( \vert V_i(B_1) \vert+ \vert V_i(B_2)  \vert \\
    && -\vert V_i(B_1\cup B_2) \vert  \big) -l_{t-1,t} \\
        & = &\sum_{i=1}^{t} k_i\vert V_i(B_1)\cap V_i(B_2)\vert -l_{t-1,t},
    \end{array}$$
 and thus the induced subgraph on $V(B_1) \cap V(B_2)$ is a tight graph. Hence, by Lemma~\ref{no_isolated}, the induced subgraph on $V(B_1) \cap V(B_2)$ has no isolated vertices, and hence, we see that $V(B_1 \cap B_2) = V(B_1)\cap V(B_2)$.

    This completes the proof of the claim.
\end{proof}
\begin{lemma}
    Let $K\in \mathbb{Z}_{>0}^{n}, L\in \mathbb{Z}_{\geq 0}^{n-1}$ be such that inequalities \eqref{conditionposet1 -k k l}, \eqref{conditionposet2 - k l l} and \eqref{conditionposet3-increasing} are valid. Let $G=(V,E)$ be a graph such that the vertex set has a partition $V=V_1\cup \dots \cup V_n$ such that all edges $e\in E$ have endpoints in $V_{i} \cup V_{i+1}$ for some $1 \leq i < n$. 
    
    If $\mathcal{M}_c$ is a matroid for all $1 \leq c\leq j$, then $\theta_{\alpha_j} = \theta_{\alpha_j'}$ for every two words $\alpha_j$ and $\alpha_j'$, such that $E_{\alpha_j}$ and $E_{\alpha_j'}$ are maximal feasible subsets supported on $V_1\cup \dots \cup V_{j+1}$. 
    \label{surplus-independent}
\end{lemma}
\begin{proof}
We are given $\alpha_j$ and $\alpha_j'$. For any word $\alpha_j = e_1, \dots e_m$ such that $E_{\alpha_j}$ is supported on $V_1\cup \dots \cup V_{j+1}$ and $d<j$, define $\alpha_{d}$ to be the word $e_1 \dots e_{\tilde{m}}$ obtained by restricting to those edges supported on $V_1\cup \dots \cup V_{d+1}$. I.e. $\alpha_{d}$ is the word resulting from deleting the edges $e_{\tilde{m}+1}\dots e_{m}$ where $\tilde{m}+1$ is the smallest index such that for all $i\geq \tilde{m}+1$, one has $e_{i} \in E(V_{d+i}\cup V_{d+i+1})$ for some $i \geq 1$. In particular, $\alpha_{0}$ will always be the empty word and for subsets $X\subseteq V_1$, that $\theta_{\alpha_0}(X) = \theta_{\alpha_0'}(X) = k_1|X|$, which serves as the basis for our induction. Suppose now that the statement holds for all $d < j$. We will show it holds for $d=j$.

Let $X\subseteq V_{j+1}$. One can verify that one can express $\theta_{\alpha_{j}}$ as a function of $\theta_{\alpha_{j-1}}$. Namely, one has:
\begin{equation*}
    \theta_{\alpha_j}(X)= \min_{\substack{W\subseteq V_{j} \cup V_{j+1}\\[2pt] X\subseteq W}}  k_{j+1}|V_{j+1} \cap W| +\theta_{\alpha_{j-1}}(V_{j}\cap W)  -|E_{\alpha_j}(W)|.
\end{equation*}

By induction, we can assume that $\theta_{\alpha_{j-1}} = \theta_{\alpha_{j-1}'}$. Thus, it suffices to consider the case where $\alpha_{j}$ and $\alpha_{j}'$ only differ for edges $e\in E(V_{j}\cup V_{j+1})$, and we may assume $\alpha_{j}\neq \alpha_{j}'$ since otherwise the lemma is trivially true. In this case, since we know that $\mathcal{M}_{j}$ is a matroid, we only need to show that in the case where $E_{\alpha_{j}} = \left (E_{\alpha_{j}'} \cup \{e\}\right) \setminus \{f\}$, one has $\theta_{\alpha_{j}}(X) = \theta_{\alpha_{j}'}(X)$, since the statement then follows for all maximally chosen words $\alpha_j$ and $\alpha_j'$. Suppose that the added edge, $e$, has endpoints $v,w$. By Lemma \ref{blocks}, we know that the intersection 
\begin{equation*}
M  = \bigcap_{\substack{T\subseteq E_{\alpha_j'},\\ v,w \in V(T),\\ T \textup{ tight }}} T
\end{equation*}
is a nonempty tight subset of $E_{\alpha_j'}$.  Furthermore, if $f$ were not in $M$, then $M$ would be a tight subset of $E_{\alpha_j}$, as $f$ is the only edge of $E_{\alpha_j'}$ not in $E_{\alpha_j}$. Then, $M \cup \{e\}$ would not be $(K,L)$-sparse, which contradicts $E_{\alpha_j}$ being $(K,L)$-sparse. It follows that $f \in M$.  Now, let $X\subseteq V$, and suppose $W$ minimises the expression for $\theta_{\alpha_j'}$, so
\begin{equation*}
\theta_{\alpha_{j}'}(X) = \sum_{i=1}^{j+1} k_i|V_i\cap W| -|E_{\alpha_j'}(W)|.
\end{equation*}

 If $f \notin E_{\alpha_j'}(W)$, then we see that

 \[|E_{\alpha_j}(W)| \geq |E_{\alpha_j'}(W)|\]

so it follows immediately that 

\begin{align*}
    \theta_{\alpha_j}(X) &\leq  \sum_{i=1}^{j+1} k_i |V_i \cap W| - |E_{\alpha_{j}}(W)| \\
    &\leq \sum_{i=1}^{j+1} k_i |V_i \cap W| - |E_{\alpha_{j}'}(W)| = \theta_{\alpha_j'}(X).
\end{align*}
 
 If $f \in E_{\alpha_j'}(W)$, then we get using Lemma \ref{Increase_when_tight - greedoids}, that
\begin{align*}
\sum_{i=1}^{j+1} k_i |V_i\cap W| - |E_{\alpha_j'}(W)|  \geq \sum_{i=1}^{j+1} k_i |V_i\cap (W\cup V(M))| - |E_{\alpha_{j}'}(W \cup V(M))|,
\end{align*}
and hence we may assume that $V(M) \subseteq W$.

Since $E_{\alpha_{j}} = \left (E_{\alpha_{j}'} \cup \{e\}\right) \setminus \{f\}$, and $e$ and $f$ are both supported on $V(M)$ it follows that,
$$|E_{\alpha_{j}'}(W\cup V(M))|=|E_{\alpha_{j}}(W\cup V(M))|.$$

It now follows that

$$\theta_{\alpha_j}(X) \leq  \sum_{i=1}^{j+1} k_i |V_i \cap W| - |E_{\alpha_{j}'}(W)| = \theta_{\alpha_j'}(X).$$

We conclude that $\theta_{\alpha_j}(X)\leq \theta_{\alpha_j'}(X)$ in all cases. Switching the roles of $E_{\alpha_j}$ and $E_{\alpha_j'}$, the same argument gives that $\theta_{\alpha_j}(X) \geq \theta_{\alpha_j'}(X)$, so $\theta_{\alpha_j'}(X) = \theta_{\alpha_j}(X)$. 
\end{proof}

\begin{lemma}\label{Matroid_1}
Let $K\in \mathbb{Z}_{>0}^{n}, L\in \mathbb{Z}_{\geq 0}^{n-1}$ be such that inequalities \eqref{conditionposet1 -k k l}, \eqref{conditionposet2 - k l l} and \eqref{conditionposet3-increasing} are valid. Let $G=(V,E)$ be a graph such that the vertex set has a partition $V=V_1\cup \dots \cup V_n$ such that all edges $e\in E$ have endpoints in $V_{i} \cup V_{i+1}$ for some $1 \leq i < n$.

Suppose that $\mathcal{M}_c$ is a matroid for all $c < j$. Then the set $\mathcal{I}_{j}$ is independent of the maximally sparse set $E^*$. Furthermore, $\mathcal{I}_j$ is non-empty whenever $E(V_{j} \cup V_{j+1})$ is non-empty.
\end{lemma}
\begin{proof}

First, assuming that $E^* \cap E(V_{j-1}\cup V_j)$ is non-empty, it is easy to see, using inequality \eqref{conditionposet2 - k l l}, that $E^* \cup \{e\}$ is $(K,L)$-sparse for any edge $e \in E(V_{j} \cup V_{j+1})$. If $E^* \cap E(V_{j-1}\cup V_j)= \emptyset$, then it is easy to see, using inequality \eqref{conditionposet1 -k k l}, that $E^* \cup \{e\}$ is $(K,L)$-sparse for any edge $e \in E(V_{j} \cup V_{j+1})$. Hence, $\mathcal{I}_j$ is non-empty whenever $E(V_{j} \cup V_{j+1})$ is non-empty.

Let $E_1^{*}$ and $E_{2}^{*}$ be two maximally sparse subsets of $E(V_1 \cup \dots \cup V_{j})$, and let $\alpha_{j-1},\alpha'_{j-1} \in \mathcal{L}_{(K,L)}$ be such that $E_{\alpha_{j-1}} = E_1^*$, $E_{\alpha'_{j-1}} = E_2^*$. Suppose that $E\subseteq E(V_{j} \cup V_{j+1})$, and that $E \cup E_1^{*}$ is sparse. Our goal is to show that $E \cup E_{2}^{*}$ is also sparse.
We make a few remarks about $E \cup E_1^{*}$ first. Let $F\subseteq E$ and consider $V_{j}(F)$ and $V_{j+1}(F)$. We note that for any subset $A \subseteq E_1^{*}$, one has, by sparsity of $F \cup A$:
\begin{equation*}
|F \cup A|=|F| +|A| \leq  \Sigma_{i=1}^{j-1}k_i |V_i(A)| + k_{j} |V_{j}(A\cup F)| +  k_{j+1} |V_{j+1}(F)| -l_{{j},{j+1}}.
\end{equation*}
Hence, one has
\begin{equation*}
|F| \leq  \Sigma_{i=1}^{j-1} k_i |V_i(A)| + k_{j} |V_{j}(A\cup F)| +  k_{j+1} |V_{j+1}(F)| -l_{{j},{j+1}} - |A|.
\end{equation*}
Since $A$ is arbitrary, we see that:
\begin{equation*}
|F| \leq   k_{j+1} |V_{j+1}(F)| -l_{j,{j+1}} + \min_{A} \{ \Sigma_{i=1}^{j-1}k_i |V_i(A)| + k_{j} |V_{j}(A\cup F)| -|A|\}.
\end{equation*}
It is then clear that: 
\begin{equation*}
|F| \leq   k_{j+1} |V_{j+1}(F)| -l_{j,{j+1}} + \theta_{\alpha_{j-1}}(V_{j}(F)).
\end{equation*}
Now, for any subset $C\subseteq E_2^{*}$, we see that
\begin{align*}
|F \cup C|&=|F| + |C| \\
&\leq  k_{j+1} |V_{j+1}(F)| -l_{j,{j+1}} + \theta_{\alpha_{j-1}}(V_{j}(C\cup F)) + |C|\\
&=   k_{j+1} |V_{j+1}(F)| -l_{j,{j+1}} + \theta_{\alpha'_{j-1}}(V_{j}(C\cup F)) + |C|\\
&\leq   k_{j+1} |V_{j+1}(F)| -l_{j,{j+1}} + \Sigma_{i=1}^{j-1} k_i|V_i(C)| + k_{j}|V_{j}(C\cup F)|\\
&- |C| + |C|\\
&\leq  \Sigma_{i=1}^{j+1} k_{i}|V_{i}(C\cup F)|-l_{j,{j+1}} 
\end{align*}
We applied the observation that if $X\subseteq Y$ then $\theta_{\alpha_{j-1}}(X) \leq \theta_{\alpha_{j-1}}(Y)$ in the second line, with $X=V_{j}(F)$ and $Y = V_{j}(C\cup F)$. We have used Lemma \ref{surplus-independent} to go from the second line to the third. Since any subset of $E_2^{*} \cup E$ is of the form $F\cup C$, this shows that $E_2^{*} \cup E$ is sparse, completing the proof of the Lemma.
\end{proof}

\begin{proof}[Proof of Theorem \ref{thm:Matroid}]

Note that $\mathcal{M}_1$ is a count matroid, so Theorem \ref{thm:Matroid} is true for $j=1$. Suppose that $\mathcal{M}_c$ is a matroid for all $c < j$. By Lemma \ref{Matroid_1}, the set $\mathcal{I}_{j}$ is independent of the maximally sparse set $E^*$. We can then fix some maximally sparse set $E^*$ such that $V(E^*) \subseteq V_1 \cup \dots \cup V_{j}$.

We will show that $\mathcal{M}_j$ is a matroid by considering its circuits. A circuit in $\mathcal{M}_j$ is a set $C \subseteq E(V_j\cup V_{j+1})$ such that there exists some subset $F^{*} \subseteq E^{*}$ with 

\[
|F^{*} \cup C| = \Sigma_{i=1}^{j+1} k_i|V_i(F^{*} \cup C)|-\ell_{j,j+1}+1
\]

but any subset $E' \subseteq E^{*} \cup C$ with $|E' \cap C| < |C|$ is $(K,L)$-sparse. 

We may assume that the collection of circuits is non-empty, since otherwise $\mathcal{M}_j$ is a uniform matroid. Note also that if $C_1$ and $C_2$ are two circuits such that $C_1 \subseteq C_2$, then $C_1 = C_2$, since any subset $E' \subseteq E^{*} \cup C_2$ with $|E' \cap C_2| < |C_2|$ is $(K,L)$-sparse. 

It remains to show that if $C_1$ and $C_2$ are two circuits such that $C_1 \neq C_2$ and $e \in C_1 \cap C_2$, then there is a circuit $C_3 \subseteq C_1 \cup C_2 \setminus \{e\}$. 

Suppose that $C_1$ and $C_2$ are two circuits, with associated subsets $F_1^{*}$ and $F_2^{*}$, such that $C_1 \neq C_2$, and let $e \in C_1 \cap C_2$. Now consider $C_1 \cup C_2 \setminus \{e\}$. Note that 
\begin{align*}
|\left((F_1^* \cup C_1) \cup (C_2 \cup F_2^*)\right) \setminus \{e\}|&=|F_1^* \cup C_1 \cup F_2^{*}  \cup C_2|-1\\
& = |F_1^* \cup C_1| + |F_2^* \cup C_2| \\
&- |\left(F_1^* \cup C_1\right) \cap \left( F_2^* \cup C_2\right)| -1.
\end{align*}
One has $\left(F_1^* \cup C_1\right) \cap \left( F_2^* \cup C_2\right) = \left(F_1^* \cap F_2^{*}\right) \cup \left( C_1 \cap C_2\right)$ as $F_{a}^{*} \cap C_{b}$ is empty for $a,b\in \{1, 2\}$. Thus one sees that the set $\left(F_1^* \cup C_1\right) \cap \left( F_2^* \cup C_2\right)$ is sparse, since $C_1\cap C_2$ is a strict subset of a circuit. By using that $C_1 \cup F_1^{*}$  and $C_2 \cup F_2^{*}$ break sparsity, one derives:
\begin{align*}
    |\left(F_1^* \cup C_1 \cup F_2^{*} \cup C_2\right) \setminus \{e\}| &\geq \sum_{i=1}^{j+1} k_i|V_i\left(\left(F_1^{*} \cup C_1 \cup F_2^{*} \cup C_2 \right) \setminus \{e\}\right) -l_{j, j+1} + 1.
\end{align*}
Thus, it follows that $C_1 \cup C_2 \setminus \{e\}$ must contain a circuit of $\mathcal{M}_j$.
\end{proof}

We are now ready to prove the main result of this article. 

\medskip
\noindent {\bf  Theorem \ref{thm:is_greedoid=true}.} {\em Let $K\in \mathbb{Z}_{>0}^{n}, L\in \mathbb{Z}_{\geq 0}^{n-1}$ be such that inequalities \eqref{conditionposet1 -k k l}, \eqref{conditionposet2 - k l l} and \eqref{conditionposet3-increasing} are valid. Let $G=(V,E)$ be a graph such that the vertex set has a partition $V=V_1\cup \dots \cup V_n$ such that all edges $e\in E$ have endpoints in $V_{i} \cup V_{i+1}$ for some $1 \leq i < n$. Then $\mathcal{L}_{(K,L)}$ is a greedoid language.}
\medskip

\begin{proof}[Proof of Theorem \ref{thm:is_greedoid=true}]
It is clear that $\mathcal{L}_{(K,L)}$ is a hereditary language. We need to show that the exchange property holds.
So suppose that $\alpha, \beta \in \mathcal{L}$, and $\vert \alpha \vert > \vert \beta \vert.$  Denote the underlying edge subset defined by $\alpha$ and $\beta$ by $E_{\alpha}$ and $E_{\beta}$ respectively.

Suppose that all edges in $E_{\beta}$ are supported on $V_1\cup\dots \cup V_j$ for some $j\geq 2$. If $E_{\beta}$ is maximally sparse on $V_1\cup \dots\cup V_j$, then any edge $e \in E_\alpha \cap E(V_{k} \cup V_{k+1})$ for the smallest $k \geq j$ with $E(V_{k} \cup V_{k+1})$ nonempty could be added to create a $(K,L)$-sparse graph, by Lemma \ref{Matroid_1}. Hence, taking any $e \in E_\alpha \cap E(V_{k} \cup V_{k+1})$ yields $\beta e \in \mathcal{L}_{(K,L)}$.

If $E_{\beta}$ is not maximally sparse on $V_1\cup \dots \cup V_j$, then note that $E_{\alpha} \cap E(V_{j-1} \cup V_j)$ and $E_{\beta} \cap E(V_{j-1} \cup V_j)$ are independent sets of the matroid $\mathcal{M}_j$. Since furthermore $|E_{\alpha} \cap E(V_{j-1} \cup V_j)| > |E_{\beta} \cap E(V_{j-1} \cup V_j)|$, there exists some $e\in E_\alpha \cap E(V_{j-1} \cap V_j)$ such that $\beta e\in \mathcal{L}_{(K,L)}$. 
\end{proof}

\section{Pebble game algorithm}\label{sec:Algorithm}
We assume that we are given a (multi)-graph $G=(V,E)$, as in Definition~\ref{KL-sparse}, and we want to test for $(K,L)$-sparsity. To do this, we introduce a pebble game algorithm.

During the pebble game, we maintain the following data
\begin{itemize}
\item A directed multigraph $D$, with $V(D)=V$, and with edge set being a subset of $E$. We call the edges of $D$ \emph{accepted} edges. For each edge $e\in D$, directed from $v$ to $w$, we will use the notation $o(e) := v$ and $t(e):=w$.
\item For each vertex $v\in V(D)$, we store a natural number $\peb(v)$, representing the number of `pebbles' on the node $v$.
\end{itemize}
The following two actions are allowed to construct the graph. They are the key steps that will be performed in the algorithm.
\begin{itemize}
\item{\AcceptEdge{$e$}: If $e=vw$ is an edge in $E(V_i \cup V_{i+1})$ that is not in $D$, and $\peb(v)+\peb(w)>l_{i,i+1}$, either
remove a pebble from $v$ (if $\peb(v) > 0$), add $(v,w)$ to $D$, and orient the edge from $v$ to $w$, or remove a pebble from $w$ (if $\peb(w) > 0$), add $(w,v)$ to $D$, and orient the edge from $w$ to $v$.} 
\item{\MovePebble{$v\leftarrow w$}: A depth-first search algorithm is done to search for a directed path whose source is $v$ and whose endpoint $w$ contains a pebble. If there exists a directed path $vf_1\dots f_kw$ of arcs $f_i$ in $D$ from $v$ to $w$, and if $\peb(w)>0$,  add a pebble to $v$, remove a pebble from $w$ and reverse the orientation for every arc $f_i$.}

\end{itemize}

We can now describe the pebble game algorithm in pseudocode in Algorithm~\ref{algo:posets}. We handle the edges in a determined order. The idea is to first process the edges between $V_1$ and $V_2$, then process the edges between $V_2$ and $V_3$, and so on. Within the sets $E(V_{i}\cup V_{i+1})$, for $i \in \{1, \dots, n-1\}$, any order is allowed.

\begin{algorithm}[h]
\textbf{Input}: A graph $G$, with a partition $V=V_1\cup \cdots \cup V_n$, such that any $e\in E$ has endpoints in $V_i$ and $V_{i+1}$ for some $1 \leq i \leq n-1$.\\
Parameters $K\in \mathbb{Z}^{n}_{>0} ,L\in \mathbb{Z}^{n-1}_{\geq 0}$, which satisfy inequalities $\eqref{conditionposet1 -k k l},\eqref{conditionposet2 - k l l}$ and $\eqref{conditionposet3-increasing}$.\\
\textbf{Output:} $(K,L)$-sparsity and/or $(K,L)$-tightness of $G$
\hrule

\begin{algorithmic} 
\State  For each vertex, set $\peb(v) \gets k_{i}$ if $v\in V_{i}$ \textit{(Initialise number of pebbles)}.
\State $D\gets(V, \emptyset)$  \textit{(Initialise $D$)}.

\For{ $i\in \{1,\cdots, n-1\}$  }
    \For{ $e=vw \in E$ with $v\in V_i$ and $w\in V_{i+1}$ }
        \While{ $e$ has not been processed }
            \If{$\peb(v)+\peb(w) > l_{i,i+1}$}
                \State \AcceptEdge{$e$}.
                \State $e$ has been processed.
            \Else
                \State Depth first search for a directed path $(v'f_1\dots f_n u)$, with $f_j\in E(D)$, $v'\in\{v,w\}$, where $u\notin \{v,w\}$, such that $\peb(u)>0$
                \If{$(v'\dots  u)$ with $\peb(u)>0$ has been found}
                    \State \MovePebble{$v'\leftarrow u$}. 
                \Else     
                    \State \Return $G$ is not sparse. 
                \EndIf
            \EndIf 
        \EndWhile
    \EndFor
\EndFor

\If{$\sum_{v\in V} \peb(v)=\max L$}
\State \Return $G$ is tight. 
\Else 
\State \Return $G$ is sparse, but not tight. 
\EndIf
\end{algorithmic}
\caption{Pebble game for $(K,L)$-sparse graphs}
\label{algo:posets}
\end{algorithm}

\newpage
Let us now prove that this pebble game recognises the $(K,L)$-sparse graphs that we have defined. For every subset $V' \subseteq V$, we define the following quantities.
\begin{align*}
    \peb(V')&:= \sum_{v\in V'} \peb(v),\\
    \spa(V')&:=\# \{e\in D \: \vert \: o(e) \in v \text{ and } t(e) \in V' \}, \\
    \out(V')&:=\# \{e\in D \: \vert \: o(e) \in v \text{ and } t(e) \notin V' \}.
\end{align*}

\begin{lemma}\label{invariants - posets}
 Let $K\in \mathbb{Z}_{>0}^{n}, L\in \mathbb{Z}_{\geq 0}^{n-1}$ be such that inequalities \eqref{conditionposet1 -k k l}, \eqref{conditionposet2 - k l l} and \eqref{conditionposet3-increasing} are valid. Let $G=(V,E)$ be a graph such that the vertex set has a partition $V=V_1\cup \dots \cup V_n$ such that all edges $e\in E$ have endpoints in $V_{i} \cup V_{i+1}$ for some $1 \leq i < n$. The following are invariants of the algorithm
    \begin{enumerate}
        \item For all $v\in V_i$, one has $\peb(v)+ \out(v) = k_{i}$.
        \item Suppose that $V'\subseteq V$. Then $$\peb(V')+ \spa(V')+ \out(V') = \sum_{i=1}^{n} k_i \vert  V' \cap V_i  \vert.$$
        \item Consider a nonempty subset $V' \subseteq V$. Then $$\peb(V') + \out(V') \geq l(V').$$ 
    \item Consider a subset $V' \subseteq V$. Then $$\spa(V') \leq \sum_{i=1}^{n} k_i \vert  V' \cap V_i \vert - l(V').$$
    \end{enumerate}
\end{lemma}
\begin{proof}
    The technique of the proof is similar to that of Lemma 10 in \cite{LEE20081425}, and we only briefly sketch the proofs of invariants $1$, $2$, and $4$ for completeness. Since invariant $3$ requires more analysis, we give the proof in more detail. 
    For the first point, this holds at the start of the algorithm. By considering what happens when moving a pebble or accepting an edge, one sees that the invariant remains valid. For invariant $2$, we see that $$\sum_{v\in V'} \peb(v) + \out(v) = \peb(V') + \spa(V') +\out(V') = \sum_{i=1}^{n} k_i \vert V' \cap V_i \vert.$$
    
     We now show that the third invariant holds after every step in the algorithm. Take a subset $V' \subseteq V$. By the convention for the maximum defining $l(V')$, if there is no $i$ with $V_{i}\neq \emptyset$ and $V_{i+1} \neq \emptyset$, then $l(V')=-\infty$, and there is nothing to show. So, we assume that $l(V')=l_{i, i+1}$ for some $i\in \{1,\cdots, n-1\}$.\\
    Suppose that an edge has been accepted with endpoints in $V'$ in the previous step of the algorithm. We need to distinguish between various cases to show the inequality. If the edge is between vertices $x,y$ with $x\in V_{i}$, and $y\in V_{i+1}$, then there must be $l_{i,i+1}$ pebbles left after accepting, since there need to have been at least $l_{i, i+1} +1$ pebbles before the edge was accepted. If the accepted edge was between vertices $x,y$ with $x\in V_{i-1}$, and $y\in V_i$, then there must be at least $l_{i-1,i}$ pebbles at $x$ and $y$, and there must be a vertex $z$ with $\peb(z) =k_{i+1}$, since the edges in $E(V_i\cup V_{i+1})$ are processed after those in $E(V_{i-1}\cup V_{i})$. Thus there must be $l_{i-1,i} + k_{i+1} > l_{i,i+1}$ pebbles left in $V'$, by inequality \eqref{conditionposet2 - k l l}. If there was an edge accepted between vertices $x,y$ with $x \in V_{k}$ and $y \in V_{k+1}$, with $k< i-1$, then there must be vertices $z$ and $w$ (by the order in which edges are accepted) with $\peb(z) =k_{i}$, $\peb(w) =k_{i+1}$, and thus there are $k_i + k_{i+1}\geq l_{i,i+1}$ pebbles left. If an edge gets accepted outside of $V'$, or if a pebble gets moved, $\peb(V') +\out(V')$ does not change, which can be seen by a case analysis.
    
    The fourth point follows immediately from the third and second invariant, by subtracting the third invariant from the second one.
\end{proof}

\begin{corollary}
 Let $K\in \mathbb{Z}_{>0}^{n}, L\in \mathbb{Z}_{\geq 0}^{n-1}$ be such that inequalities \eqref{conditionposet1 -k k l}, \eqref{conditionposet2 - k l l} and \eqref{conditionposet3-increasing} are valid. Let $G=(V,E)$ be a graph such that the vertex set has a partition $V=V_1\cup \dots \cup V_n$ such that all edges $e\in E$ have endpoints in $V_{i} \cup V_{i+1}$ for some $1 \leq i < n$.  
 
 If all edges of $G$ are accepted in the pebble game, then $G$ is  $(K,L)$-sparse. If there are exactly $l_{n-1,n}$ pebbles left in the graph, then $G$ is tight.
\end{corollary}
\begin{proof}
This is an easy consequence of the fourth and second invariants of Lemma \ref{invariants - posets}.
\end{proof}

We have now shown that if every edge in $G$ is accepted in the pebble game, then $G$ is sparse. Now we show the converse, namely that if $G$ is sparse, then every edge of $G$ will be accepted.\\ In the following proof, we need to define the notion of the $Reach$ of a vertex in an oriented graph $D$. Let 
\[
\textup{Reach}(v):=\{v\} \cup \{w\in V(D) \: \vert \: \mbox{ there exists an oriented path from } v \mbox{ to } w\}.
\]

\begin{lemma}\label{lemma-add-edge(poset)}
 Let $K\in \mathbb{Z}_{>0}^{n}, L\in \mathbb{Z}_{\geq 0}^{n-1}$ be such that inequalities \eqref{conditionposet1 -k k l}, \eqref{conditionposet2 - k l l} and \eqref{conditionposet3-increasing} are valid. Let $G=(V,E)$ be a graph such that the vertex set has a partition $V=V_1\cup \dots \cup V_n$ such that all edges $e\in E$ have endpoints in $V_{i} \cup V_{i+1}$ for some $1 \leq i < n$. 
 
Suppose that after running some steps of the algorithm $D$ is the set of edges we accepted, where all edges in $D$ are between $V_a$ and $V_{a+1}$ for $a\leq i$. Suppose that $e$ is an edge of $G$ which has not been accepted yet such that $e=(uv)$, where $u\in V_i$ and $v\in V_{i+1}$. Then the edge can be accepted (after perhaps moving some pebbles).
\end{lemma}
\begin{proof}
If there are more than $l_{i, i+1}$ pebbles at $u$ and $v$ together, we can add the edge and remove a pebble. So, we may assume that $\peb(u)+\peb(v)<l_{i, i+1}+1$.
Let $V'=\textup{Reach}(u)\cup \textup{Reach}(v)$, and let $E_G(V')$ denote the edges induced by $V'$ in $G$. Note that $V'\cap V_i$ and $V'\cap  V_{i+1}$ are non-empty since $u\in V'\cap V_i$ and $v\in V'\cap V_{i+1}$. Furthermore $V'\cap V_k =\emptyset$, for $k\geq i+2$, since these edges get processed at some later point. Since $G$ is sparse, $$\spa(V')< \vert E_G(V') \vert \leq \sum_{j=1}^{i+1}k_j\vert V' \cap V_j  \vert - l_{i,i+1},$$ where the first inequality is strict since $e\in  E_G(V') \setminus \spa(V')$, and $\max\{l_{k,k+1}\}=l_{i,i+1}$ since the $l_{k,k+1}$ are increasing. By the definition of $\textup{Reach}$, $\out(V')=0$ and by invariant 2 in Lemma~\ref{invariants - posets}, 
\[
\spa(V')+ \peb(V')= \sum_{j=1}^{i+1}k_j\vert V_j \cap V' \vert
\]
and thus
\[
\peb(V') > l_{i, i+1} .
\]
Hence, there is a vertex $w\in V'$ such that $w\notin \{v,u\}$ with at least one pebble left. Since it is reachable by either $v$ or $u$, we can move that pebble (changing the direction of the path). Repeat this process until we get enough pebbles together at $u$ and $v$ so that we can accept the edge $(uv)$ and add it to $D$.
\end{proof}

\begin{corollary}
  Let $K\in \mathbb{Z}_{>0}^{n}, L\in \mathbb{Z}_{\geq 0}^{n-1}$ be such that inequalities \eqref{conditionposet1 -k k l}, \eqref{conditionposet2 - k l l} and \eqref{conditionposet3-increasing} are valid. Let $G=(V,E)$ be a graph such that the vertex set has a partition $V=V_1\cup \dots \cup V_n$ such that all edges $e\in E$ have endpoints in $V_{i} \cup V_{i+1}$ for some $1 \leq i < n$. 
  
  If the graph $G$ is $(K,L)$-sparse, then the pebble game accepts all edges of $G$.
\end{corollary}
\begin{proof}
This follows easily from Lemma~\ref{lemma-add-edge(poset)}.
\end{proof}

Note that if $l_{i,i+1}$ is the same for all parts $V_i$, then the structure is a count matroid, and the pebble game recognises the independent sets of this count matroid. In particular, it can recognise $(k_1,k_2,l)$-sparse bipartite graphs. In this way, our algorithm extends Berg and Jord\'an's algorithm for the $k$-plane matroid. Running the algorithm on a graph that is not sparse finds the maximal feasible sets introduced in Section \ref{sec:greedoid}, and not all maximally sparse subsets.

Finally, let us consider the complexity of the algorithm. There are at most $l_{n-1, n}$ searches for each edge, and depth first search takes $\mathcal{O}(|V|+|E|)$ time, which is performed at most $|E|$ times. Hence, one has a time complexity of 
\begin{equation*}
    \mathcal{O}\left(l_{n-1, n}|E|\left(|V|+|E|\right)\right).
\end{equation*}
For the space complexity, the graph $G$ is stored. For each vertex, we store the number of pebbles, and for the edges, we store their direction. Hence, the space complexity is given by
\begin{equation*}
\mathcal{O}(|V|+|E|).
\end{equation*}

\section{Conclusions and open problems}\label{sec:conc,prob}


    We conclude by pointing towards avenues for further work. In this article, we have proved that the pebble game recognises greedoid structures. However, one could consider the same algorithm for graphs where one first processes a subset of edges $E_1$, then $E_2$, and so on, using different $\ell$ for processing the edges. A natural question is whether this always defines a greedoid and whether this greedoid can be described combinatorially, as we have in this paper.
    
One classical feature of sparsity matroids is that they imply the existence of graph decompositions, as illustrated in the Nash-Williams theorem. Our Theorem \ref{arboricity} also provides an example of this. A natural question from this perspective is whether one can find a similar result for $(K,L)$-sparsity. It seems natural to study arboricity in this setting, i.e., decompositions into rooted trees, where one demands that the root of the tree belongs to $V_1$. 

Concretely, in the case where $K=(1, \dots, 1)$ and $L=(1, \dots, 1)$, one easily sees that maximally sparse graphs are unions of rooted trees with roots in $V_1$, such that any tree remains a tree when looking at the edges supported on $V_1\cup \dots \cup V_j$ for each $j\in \{2, \dots, n\}$. One could then ask whether maximally $(K,L)$-sparse graphs with $K=(k, \dots, k)$ and $L=(k, \dots, k)$ have a decomposition into $k$ of these forests.
    
Finally, one feature of sparsity matroids, useful in rigidity theory, is that one can describe how sparse graphs are inductively constructed for certain values of $k$ and $\ell$ \cite{SZEGO20061211}. While we have not investigated this in this work, an avenue for research is whether one can inductively construct the sparse graphs for certain values of $K$ and $L$.

\section*{Acknowledgements}
The authors are grateful to Tibor Jord\'an for pointing to reference \cite{FRANK}. 

This work has been supported by the Knut and Alice Wallenberg Foundation Grant 2020.0001 and 2020.0007, and Kempe foundation no. SMK21-0074. This work was partially supported by the Wallenberg AI, Autonomous Systems and
Software Program (WASP) funded by the Knut and Alice Wallenberg Foundation. Signe Lundqvist was partially supported by the FWO grants G0F5921N (Odysseus) and G023721N, and by the KU Leuven grant iBOF/23/064.

\bibliographystyle{abbrv}
\bibliography{Bibliography}

@InProceedings{Berg-Jordan,
author="Berg, Alex R.
and Jord{\'a}n, Tibor",
editor="Di Battista, Giuseppe
and Zwick, Uri",
title="Algorithms for Graph Rigidity and Scene Analysis",
booktitle="Algorithms - ESA 2003",
year="2003",
publisher="Springer Berlin Heidelberg",
address="Berlin, Heidelberg",
pages="78--89",
isbn="978-3-540-39658-1"
}

@article{Frank2012,
  title = {Connections in Combinatorial Optimization},
  volume = {160},
  ISSN = {0166-218X},
  url = {http://dx.doi.org/10.1016/j.dam.2011.09.003},
  DOI = {10.1016/j.dam.2011.09.003},
  number = {12},
  journal = {Discrete Applied Mathematics},
  publisher = {Elsevier BV},
  author = {Frank,  András},
  year = {2012},
  pages = {1875}
}

@article{STREINU20091944,
title = {Sparse hypergraphs and pebble game algorithms},
journal = {European Journal of Combinatorics},
volume = {30},
number = {8},
pages = {1944-1964},
year = {2009},
issn = {0195-6698},
doi = {https://doi.org/10.1016/j.ejc.2008.12.018},
url = {https://www.sciencedirect.com/science/article/pii/S0195669808002795},
author = {Ileana Streinu and Louis Theran}
}

@article{LEE20081425,
title = {Pebble game algorithms and sparse graphs},
journal = {Discrete Mathematics},
volume = {308},
number = {8},
pages = {1425-1437},
year = {2008},
note = {Third European Conference on Combinatorics},
issn = {0012-365X},
doi = {https://doi.org/10.1016/j.disc.2007.07.104},
url = {https://www.sciencedirect.com/science/article/pii/S0012365X07005602},
author = {Audrey Lee and Ileana Streinu}
}

@article {Laman,
    AUTHOR = {Laman, G.},
     TITLE = {On graphs and rigidity of plane skeletal structures},
   JOURNAL = {J. Engrg. Math.},
  FJOURNAL = {Journal of Engineering Mathematics},
    VOLUME = {4},
      YEAR = {1970},
     PAGES = {331--340},
      ISSN = {0022-0833},
   MRCLASS = {05.40},
  MRNUMBER = {269535},
       DOI = {10.1007/BF01534980},
       URL = {https://doi.org/10.1007/BF01534980},
}

@article{Geiringer,
  title={{\"U}ber die Gliederung ebener Fachwerke},
  author={Hilda Pollaczek-Geiringer},
  journal={Zamm-zeitschrift Fur Angewandte Mathematik Und Mechanik},
  year={1927},
  volume={7},
  pages={58-72}
}

@article {Pebble_original,
    AUTHOR = {Jacobs, Donald J. and Hendrickson, Bruce},
     TITLE = {An algorithm for two-dimensional rigidity percolation: the
              pebble game},
   JOURNAL = {J. Comput. Phys.},
  FJOURNAL = {Journal of Computational Physics},
    VOLUME = {137},
      YEAR = {1997},
    NUMBER = {2},
     PAGES = {346--365},
      ISSN = {0021-9991},
   MRCLASS = {65Y25 (82C99)},
  MRNUMBER = {1481894},
       DOI = {10.1006/jcph.1997.5809},
       URL = {https://doi.org/10.1006/jcph.1997.5809},
}

@article {Plane_Matroid,
    AUTHOR = {Whiteley, Walter},
     TITLE = {A matroid on hypergraphs, with applications in scene analysis
              and geometry},
   JOURNAL = {Discrete Comput. Geom.},
  FJOURNAL = {Discrete \& Computational Geometry. An International Journal
              of Mathematics and Computer Science},
    VOLUME = {4},
      YEAR = {1989},
    NUMBER = {1},
     PAGES = {75--95},
      ISSN = {0179-5376},
   MRCLASS = {05B35 (05C65 52A37)},
  MRNUMBER = {964145},
MRREVIEWER = {Robert Connelly},
       DOI = {10.1007/BF02187716},
       URL = {https://doi.org/10.1007/BF02187716},
}

@incollection {discrete_matroids,
    AUTHOR = {Whiteley, Walter},
     TITLE = {Some matroids from discrete applied geometry},
 BOOKTITLE = {Matroid theory ({S}eattle, {WA}, 1995)},
    SERIES = {Contemp. Math.},
    VOLUME = {197},
     PAGES = {171--311},
 PUBLISHER = {Amer. Math. Soc., Providence, RI},
      YEAR = {1996},
   MRCLASS = {05B35 (52B40 55U15 68U10)},
  MRNUMBER = {1411692},
MRREVIEWER = {Tiong Seng Tay},
       DOI = {10.1090/conm/197/02540},
       URL = {https://doi.org/10.1090/conm/197/02540},
}

@article {Gabow-Westerman,
    AUTHOR = {Gabow, Harold N. and Westermann, Herbert H.},
     TITLE = {Forests, frames, and games: algorithms for matroid sums and
              applications},
   JOURNAL = {Algorithmica},
  FJOURNAL = {Algorithmica. An International Journal in Computer Science},
    VOLUME = {7},
      YEAR = {1992},
    NUMBER = {5-6},
     PAGES = {465--497},
      ISSN = {0178-4617},
   MRCLASS = {90C27 (05B35 68Q25 68R10 90D46)},
  MRNUMBER = {1154585},
MRREVIEWER = {Ulrich Faigle},
       DOI = {10.1007/BF01758774},
       URL = {https://doi.org/10.1007/BF01758774},
}

@article{madarasi2025efficient,
  title={Efficient Algorithms and Implementations for Extracting Maximum-Size $(k,\ell) $-Sparse Subgraphs},
  author={Madarasi, P{\'e}ter},
  journal={arXiv preprint arXiv:2511.16877},
  year={2025}
}

@article{IMAI198579,
title = {On combinatorial structures of line drawings of polyhedra},
journal = {Discrete Applied Mathematics},
volume = {10},
number = {1},
pages = {79-92},
year = {1985},
issn = {0166-218X},
doi = {https://doi.org/10.1016/0166-218X(85)90060-5},
url = {https://www.sciencedirect.com/science/article/pii/0166218X85900605},
author = {Hiroshi Imai}
}

@article{SZEGO20061211,
title = {On constructive characterizations of (k,l)-sparse graphs},
journal = {European Journal of Combinatorics},
volume = {27},
number = {7},
pages = {1211-1223},
year = {2006},
issn = {0195-6698},
doi = {https://doi.org/10.1016/j.ejc.2006.06.016},
url = {https://www.sciencedirect.com/science/article/pii/S0195669806001193},
author = {László Szegő}
}

@article{Streinu2010,
  title = {Slider-Pinning Rigidity: a Maxwell–Laman-Type Theorem},
  volume = {44},
  ISSN = {1432-0444},
  url = {http://dx.doi.org/10.1007/s00454-010-9283-y},
  DOI = {10.1007/s00454-010-9283-y},
  number = {4},
  journal = {Discrete \& Computational Geometry},
  publisher = {Springer Science and Business Media LLC},
  author = {Streinu,  Ileana and Theran,  Louis},
  year = {2010},
  pages = {812–837}
}

@article{lee2007graded,
  title={Graded Sparse Graphs and Matroids},
  author={Lee, Audrey and Streinu, Ileana and Theran, Louis},
  journal={Journal of Universal Computer Science},
  volume={13},
  number={11},
  pages={1671--1679},
  year={2007}
}

@article{LOREA1979103,
title = {On matroidal families},
journal = {Discrete Mathematics},
volume = {28},
number = {1},
pages = {103-106},
year = {1979},
issn = {0012-365X},
doi = {https://doi.org/10.1016/0012-365X(79)90190-0},
url = {https://www.sciencedirect.com/science/article/pii/0012365X79901900},
author = {M. Lorea}
}

@incollection {k-planematroids,
    AUTHOR = {Servatius, Brigitte},
     TITLE = {{$k$}-plane matroids and {W}hiteley's flattening conjectures},
 BOOKTITLE = {Combinatorics, graph theory and computing},
    SERIES = {Springer Proc. Math. Stat.},
    VOLUME = {388},
     PAGES = {109--115},
 PUBLISHER = {Springer, Cham},
      YEAR = {2022},
      ISBN = {978-3-031-05374-0; 978-3-031-05375-7},
   MRCLASS = {05B35 (52B40 52C25)},
  MRNUMBER = {4508519},
       DOI = {10.1007/978-3-031-05375-7\_7},
       URL = {https://doi.org/10.1007/978-3-031-05375-7_7},
}

@article{FRANK,
title = {On decomposing a hypergraph into k connected sub-hypergraphs},
journal = {Discrete Applied Mathematics},
volume = {131},
number = {2},
pages = {373-383},
year = {2003},
issn = {0166-218X},
doi = {https://doi.org/10.1016/S0166-218X(02)00463-8},
url = {https://www.sciencedirect.com/science/article/pii/S0166218X02004638},
author = {András Frank and Tamás Király and Matthias Kriesell},
}

@article {Nash-Williams,
    AUTHOR = {Nash-Williams, C. St. J. A.},
     TITLE = {Edge-disjoint spanning trees of finite graphs},
   JOURNAL = {J. London Math. Soc.},
  FJOURNAL = {The Journal of the London Mathematical Society},
    VOLUME = {36},
      YEAR = {1961},
     PAGES = {445--450},
      ISSN = {0024-6107,1469-7750},
   MRCLASS = {05.45},
  MRNUMBER = {133253},
MRREVIEWER = {W.\ T.\ Tutte},
       DOI = {10.1112/jlms/s1-36.1.445},
       URL = {https://doi.org/10.1112/jlms/s1-36.1.445},
}

@InProceedings{Greedoids_1,
author="Korte, B.
and Lov{\'a}sz, L.",
title="Mathematical structures underlying greedy algorithms",
booktitle="Fundamentals of Computation Theory",
year="1981",
publisher="Springer Berlin Heidelberg",
pages="205--209",
isbn="978-3-540-38765-7"
}

@book{Greedoids_2,
    AUTHOR = {Korte, Bernhard and Lov\'{a}sz, L\'{a}szl\'{o} and Schrader,
              Rainer},
     TITLE = {Greedoids},
    SERIES = {Algorithms and Combinatorics},
    VOLUME = {4},
 PUBLISHER = {Springer-Verlag, Berlin},
      YEAR = {1991},
     PAGES = {viii+211},
      ISBN = {3-540-18190-3},
   MRCLASS = {90-02 (05B35 90C10 90C27)},
  MRNUMBER = {1183735},
MRREVIEWER = {Ulrich\ Faigle},
       DOI = {10.1007/978-3-642-58191-5},
       URL = {https://doi.org/10.1007/978-3-642-58191-5},
}

\end{document}